\documentstyle[11pt]{article}
\begin{document}
%
%
%
\newtheorem{theorem}      {Th\'eor\`eme}[section]
\newtheorem{theorem*}     {theorem}
\newtheorem{proposition}  [theorem]{Proposition}
\newtheorem{definition}   [theorem]{Definition}
\newtheorem{e-lemme}        [theorem]{Lemma}
\newtheorem{cor}   [theorem]{Corollaire}
\newtheorem{resultat}     [theorem]{R\'esultat}
\newtheorem{eexercice}    [theorem]{Exercice}
\newtheorem{rrem}    [theorem]{Remarque}
\newtheorem{pprobleme}    [theorem]{Probl\`eme}
\newtheorem{eexemple}     [theorem]{Exemple}
\newcommand{\preuve}      {\paragraph{Preuve}}
\newenvironment{probleme} {\begin{pprobleme}\rm}{\end{pprobleme}}
\newenvironment{remarque} {\begin{rremarque}\rm}{\end{rremarque}}
\newenvironment{exercice} {\begin{eexercice}\rm}{\end{eexercice}}
\newenvironment{exemple}  {\begin{eexemple}\rm}{\end{eexemple}}
%
%
\newtheorem{e-theo}      [theorem]{Theorem}
\newtheorem{theo*}     [theorem]{Theorem}
\newtheorem{e-pro}  [theorem]{Proposition}
\newtheorem{e-def}   [theorem]{Definition}
\newtheorem{e-lem}        [theorem]{Lemma}
\newtheorem{e-cor}   [theorem]{Corollary}
\newtheorem{e-resultat}     [theorem]{Result}
\newtheorem{ex}    [theorem]{Exercise}
\newtheorem{e-rem}    [theorem]{Remark}
\newtheorem{prob}    [theorem]{Problem}
\newtheorem{example}     [theorem]{Example}
\newcommand{\proof}         {\paragraph{Proof~: }}
\newcommand{\hint}          {\paragraph{Hint}}
\newcommand{\heuristicproof}{\paragraph{heuristic proof}}
\newenvironment{e-probleme} {\begin{e-pprobleme}\rm}{\end{e-pprobleme}}
\newenvironment{e-remarque} {\begin{e-rremarque}\rm}{\end{e-rremarque}}
\newenvironment{e-exercice} {\begin{e-eexercice}\rm}{\end{e-eexercice}}
\newenvironment{e-exemple}  {\begin{e-eexemple}\rm}{\end{e-eexemple}}
\newcommand{\reell}    {{{\rm I\! R}^l}}
\newcommand{\reeln}    {{{\rm I\! R}^n}}
\newcommand{\reelk}    {{{\rm I\! R}^k}}
\newcommand{\reelm}    {{{\rm I\! R}^m}}
\newcommand{\reelp}    {{{\rm I\! R}^p}}
\newcommand{\reeld}    {{{\rm I\! R}^d}}
\newcommand{\reeldd}   {{{\rm I\! R}^{d\times d}}}
\newcommand{\reelnn}   {{{\rm I\! R}^{n\times n}}}
\newcommand{\reelnd}   {{{\rm I\! R}^{n\times d}}}
\newcommand{\reeldn}   {{{\rm I\! R}^{d\times n}}}
\newcommand{\reelkd}   {{{\rm I\! R}^{k\times d}}}
\newcommand{\reelkl}   {{{\rm I\! R}^{k\times l}}}
\newcommand{\reelN}    {{{\rm I\! R}^N}}
\newcommand{\reelM}    {{{\rm I\! R}^M}}
\newcommand{\reelplus} {{{\rm I\! R}^+}}
\newcommand{\reelo}    {{{\rm I\! R}\setminus\{0\}}}
\newcommand{\reld}    {{{\rm I\! R}_d}}
\newcommand{\relplus} {{{\rm I\! R}_+}}
\newcommand{\1}        {{\bf 1}}

\newcommand{\cov}      {{\hbox{cov}}}
\newcommand{\sss}      {{\cal S}}
\newcommand{\indic}    {{{\rm I\!\! I}}}
\newcommand{\pp}       {{{\rm I\!\!\! P}}}
\newcommand{\qq}       {{{\rm I\!\!\! Q}}}
\newcommand{\ee}       {{{\rm I\! E}}}

\newcommand{\B}        {{{\rm I\! B}}}
\newcommand{\cc}       {{{\rm I\!\!\! C}}}
\newcommand{\HHH}        {{{\rm I\! H}}}
\newcommand{\N}        {{{\rm I\! N}}}
\newcommand{\R}        {{{\rm I\! R}}}
\newcommand{\D}        {{{\rm I\! D}}}
\newcommand{\Z}       {{{\rm Z\!\! Z}}}
\newcommand{\C}        {{\bf C}}	
\newcommand{\T}        {{\bf T}}	
\newcommand{\E}        {{\bf E}}	
\newcommand{\rfr}[1]    {\stackrel{\circ}{#1}}
\newcommand{\equiva}    {\displaystyle\mathop{\simeq}}
\newcommand{\eqdef}     {\stackrel{\triangle}{=}}
\newcommand{\limps}     {\mathop{\hbox{\rm lim--p.s.}}}
\newcommand{\Limsup}    {\mathop{\overline{\rm lim}}}
\newcommand{\Liminf}    {\mathop{\underline{\rm lim}}}
\newcommand{\Inf}       {\mathop{\rm Inf}}
\newcommand{\vers}      {\mathop{\;{\rightarrow}\;}}
\newcommand{\versup}    {\mathop{\;{\nearrow}\;}}
\newcommand{\versdown}  {\mathop{\;{\searrow}\;}}
\newcommand{\vvers}     {\mathop{\;{\longrightarrow}\;}}
\newcommand{\cvetroite} {\mathop{\;{\Longrightarrow}\;}}
\newcommand{\ieme}      {\hbox{i}^{\hbox{\smalltype\`eme}}}
\newcommand{\eqps}      {\, \buildrel \rm \hbox{\rm\smalltype p.s.} \over = \,}
\newcommand{\eqas}      {\,\buildrel\rm\hbox{\rm\smalltype a.s.} \over = \,}
\newcommand{\argmax}    {\hbox{{\rm Arg}}\max}
\newcommand{\argmin}    {\hbox{{\rm Arg}}\min}
\newcommand{\indep}{\perp\!\!\!\!\perp}
\newcommand{\abs}[1]{\left| #1 \right|}
\newcommand{\crochet}[2]{\langle #1 \,,\, #2 \rangle}
\newcommand{\espc}[3]   {E_{#1}\left(\left. #2 \right| #3 \right)}
\newcommand{\rang}{\hbox{rang}}
\newcommand{\rank}{\hbox{rank}}
\newcommand{\signe}{\hbox{signe}}
\newcommand{\sign}{\hbox{sign}}

\newcommand\hA{{\widehat A}}
\newcommand\hB{{\widehat B}}
\newcommand\hC{{\widehat C}}
\newcommand\hD{{\widehat D}}
\newcommand\hE{{\widehat E}}
\newcommand\hF{{\widehat F}}
\newcommand\hG{{\widehat G}}
\newcommand\hH{{\widehat H}}
\newcommand\hI{{\widehat I}}
\newcommand\hJ{{\widehat J}}
\newcommand\hK{{\widehat K}}
\newcommand\hL{{\widehat L}}
\newcommand\hM{{\widehat M}}
\newcommand\hN{{\widehat N}}
\newcommand\hO{{\widehat O}}
\newcommand\hP{{\widehat P}}
\newcommand\hQ{{\widehat Q}}
\newcommand\hR{{\widehat R}}
\newcommand\hS{{\widehat S}}
\newcommand\hTT{{\widehat T}}
\newcommand\hU{{\widehat U}}
\newcommand\hV{{\widehat V}}
\newcommand\hW{{\widehat W}}
\newcommand\hX{{\widehat X}}
\newcommand\hY{{\widehat Y}}
\newcommand\hZ{{\widehat Z}}

\newcommand\ha{{\widehat a}}
\newcommand\hb{{\widehat b}}
\newcommand\hc{{\widehat c}}
\newcommand\hd{{\widehat d}}
\newcommand\he{{\widehat e}}
\newcommand\hf{{\widehat f}}
\newcommand\hg{{\widehat g}}
\newcommand\hh{{\widehat h}}
\newcommand\hi{{\widehat i}}
\newcommand\hj{{\widehat j}}
\newcommand\hk{{\widehat k}}
\newcommand\hl{{\widehat l}}
\newcommand\hm{{\widehat m}}
\newcommand\hn{{\widehat n}}
\newcommand\ho{{\widehat o}}
\newcommand\hp{{\widehat p}}
\newcommand\hq{{\widehat q}}
\newcommand\hr{{\widehat r}}
\newcommand\hs{{\widehat s}}
\newcommand\htt{{\widehat t}}
\newcommand\hu{{\widehat u}}
\newcommand\hv{{\widehat v}}
\newcommand\hw{{\widehat w}}
\newcommand\hx{{\widehat x}}
\newcommand\hy{{\widehat y}}
\newcommand\hz{{\widehat z}}

\newcommand\tA{{\widetilde A}}
\newcommand\tB{{\widetilde B}}
\newcommand\tC{{\widetilde C}}
\newcommand\tD{{\widetilde D}}
\newcommand\tE{{\widetilde E}}
\newcommand\tF{{\widetilde F}}
\newcommand\tG{{\widetilde G}}
\newcommand\tH{{\widetilde H}}
\newcommand\tI{{\widetilde I}}
\newcommand\tJ{{\widetilde J}}
\newcommand\tK{{\widetilde K}}
\newcommand\tL{{\widetilde L}}
\newcommand\tM{{\widetilde M}}
\newcommand\tN{{\widetilde N}}
\newcommand\tOO{{\widetilde O}}
\newcommand\tP{{\widetilde P}}
\newcommand\tQ{{\widetilde Q}}
\newcommand\tR{{\widetilde R}}
\newcommand\tS{{\widetilde S}}
\newcommand\tTT{{\widetilde T}}
\newcommand\tU{{\widetilde U}}
\newcommand\tV{{\widetilde V}}
\newcommand\tW{{\widetilde W}}
\newcommand\tX{{\widetilde X}}
\newcommand\tY{{\widetilde Y}}
\newcommand\tZ{{\widetilde Z}}

\newcommand\ta{{\widetilde a}}
\newcommand\tb{{\widetilde b}}
\newcommand\tc{{\widetilde c}}
\newcommand\td{{\widetilde d}}
\newcommand\te{{\widetilde e}}
\newcommand\tf{{\widetilde f}}
\newcommand\tg{{\widetilde g}}
\newcommand\th{{\widetilde h}}
\newcommand\ti{{\widetilde i}}
\newcommand\tj{{\widetilde j}}
\newcommand\tk{{\widetilde k}}
\newcommand\tl{{\widetilde l}}
\newcommand\tm{{\widetilde m}}
\newcommand\tn{{\widetilde n}}
\newcommand\tio{{\widetilde o}}
\newcommand\tp{{\widetilde p}}
\newcommand\tq{{\widetilde q}}
\newcommand\tr{{\widetilde r}}
\newcommand\ts{{\widetilde s}}
\newcommand\tit{{\widetilde t}}
\newcommand\tu{{\widetilde u}}
\newcommand\tv{{\widetilde v}}
\newcommand\tw{{\widetilde w}}
\newcommand\tx{{\widetilde x}}
\newcommand\ty{{\widetilde y}}
\newcommand\tz{{\widetilde z}}

\newcommand\bA{{\overline A}}
\newcommand\bB{{\overline B}}
\newcommand\bC{{\overline C}}
\newcommand\bD{{\overline D}}
\newcommand\bE{{\overline E}}
\newcommand\bFF{{\overline F}}
\newcommand\bG{{\overline G}}
\newcommand\bH{{\overline H}}
\newcommand\bI{{\overline I}}
\newcommand\bJ{{\overline J}}
\newcommand\bK{{\overline K}}
\newcommand\bL{{\overline L}}
\newcommand\bM{{\overline M}}
\newcommand\bN{{\overline N}}
\newcommand\bO{{\overline O}}
\newcommand\bP{{\overline P}}
\newcommand\bQ{{\overline Q}}
\newcommand\bR{{\overline R}}
\newcommand\bS{{\overline S}}
\newcommand\bT{{\overline T}}
\newcommand\bU{{\overline U}}
\newcommand\bV{{\overline V}}
\newcommand\bW{{\overline W}}
\newcommand\bX{{\overline X}}
\newcommand\bY{{\overline Y}}
\newcommand\bZ{{\overline Z}}

\newcommand\ba{{\overline a}}
\newcommand\bb{{\overline b}}
\newcommand\bc{{\overline c}}
\newcommand\bd{{\overline d}}
\newcommand\be{{\overline e}}
\newcommand\bff{{\overline f}}
\newcommand\bg{{\overline g}}
\newcommand\bh{{\overline h}}
\newcommand\bi{{\overline i}}
\newcommand\bj{{\overline j}}
\newcommand\bk{{\overline k}}
\newcommand\bl{{\overline l}}
\newcommand\bm{{\overline m}}
\newcommand\bn{{\overline n}}
\newcommand\bo{{\overline o}}
\newcommand\bp{{\overline p}}
\newcommand\bq{{\overline q}}
\newcommand\br{{\overline r}}
\newcommand\bs{{\overline s}}
\newcommand\bt{{\overline t}}
\newcommand\bu{{\overline u}}
\newcommand\bv{{\overline v}}
\newcommand\bw{{\overline w}}
\newcommand\bx{{\overline x}}
\newcommand\by{{\overline y}}
\newcommand\bz{{\overline z}}

%
\newcommand{\AAA}{{\cal A}}
\newcommand{\BB}{{\cal B}}
\newcommand{\CC}{{\cal C}}
\newcommand{\DD}{{\cal D}}
\newcommand{\EE}{{\cal E}}
\newcommand{\FF}{{\cal F}}
\newcommand{\GG}{{\cal G}}
\newcommand{\HH}{{\cal H}}
\newcommand{\II}{{\cal I}}
\newcommand{\JJ}{{\cal J}}
\newcommand{\KK}{{\cal K}}
\newcommand{\LL}{{\cal L}}
\newcommand{\NN}{{\cal N}}
\newcommand{\MM}{{\cal M}}
\newcommand{\OO}{{\cal O}}
\newcommand{\PP}{{\cal P}}
\newcommand{\QQ}{{\cal Q}}
\newcommand{\RR}{{\cal R}}
\newcommand{\SS}{{\cal S}}
\newcommand{\TT}{{\cal T}}
\newcommand{\UU}{{\cal U}}
\newcommand{\VV}{{\cal V}}
\newcommand{\WW}{{\cal W}}
\newcommand{\XX}{{\cal X}}
\newcommand{\YY}{{\cal Y}}
\newcommand{\ZZ}{{\cal Z}}
\newcommand{\tbullet}{$\bullet$}
\newcommand{\ot}{\leftarrow}
\newcommand{\newblock}{}
\newcommand{\carre}{\hfill$\Box$}
\newcommand{\carreb}{\hfill\rule{0.25cm}{0.25cm}}
%
%
\newcommand{\dontforget}[1]
{{\mbox{}\\\noindent\rule{1cm}{2mm}\hfill don't forget : #1 \hfill\rule{1cm}{2mm}}\typeout{---------- don't forget : #1 ------------}}
\newcommand{\note}[2]
{ \noindent{\sf #1 \hfill \today}

\noindent\mbox{}\hrulefill\mbox{}
\begin{quote}\begin{quote}\sf #2\end{quote}\end{quote}
\noindent\mbox{}\hrulefill\mbox{}
\vspace{1cm}
}
\newcommand{\rond}[1]     {\stackrel{\circ}{#1}}
\newcommand{\rondf}       {\stackrel{\circ}{\FF}}
\newcommand{\point}[1]     {\stackrel{\cdot}{#1}}

\newcommand\relatif{{\rm \rlap Z\kern 3pt Z}}

\def\diagram#1{\def\normalbaselines{\baselineskip=Opt
\lineskip=10Opt\lineskiplimit=1pt}  \matrix{#1}}

\title{\huge  On partial analyticity  of CR mappings  }
\author{}
\date{}
\maketitle

\begin{center}
Bernard COUPET, Sergey PINCHUK\footnote{ The second author is partially supported 
by the NSF grant DMS 96 225 94} and Alexandre SUKHOV
\end{center}
\bigskip

{\small
Abstract. We study the problem of holomorphic extension of a smooth CR mapping from a
real analytic hypersurface to a real algebraic set in complex spaces of different
dimensions.

AMS Mathematics Subject Classification: 32D15,32D99,32F25,32H99 

Key words: CR mapping, holomorphic extension, reflection principle,  transcendence
degree }

\section{Introduction}

 Let $X$ and $Y$ be  real analytic Cauchy - Riemann
manifolds  in complex affine spaces (of different dimensions, in general), 
and $f: X \longrightarrow Y$ be a smooth ${\cal CR}$ mapping. It is natural to ask
under what conditions $f$ is real analytic (and, therefore, extends
holomorphically to a neighborhood of $X$)?  The intrinsic tool to study this problem 
is  the reflection principle. It has two major variations: the analytic one (which
makes use of the  tangent Cauchy - Riemann operators) and the geometric one
(involving a study of analytic geometry of the Segre families). In  the present paper
we assume that $Y$ is a real {\it algebraic} variety and use the analytic approach.
Algebraicity of $Y$ allows to involve additionaly methods  of commutative algebra. A
local character of the problem makes it more convenient to formulate some of our
results in terms of germs of ${\cal CR}$ mappings.

Let $X$ be a real analytic  hypersurface in $\cc^n$, minimal at a point $p \in X$. As
usual, the minimality means that $X$ contains no germs of complex hypersurfaces in a
neighborhood of $p$. Let also $Y$ be a real algebraic subset of $\cc^{N}$ i.e. the
zero set of finite number of real valued polynomials in $\cc^N$. Fix a point $p \in
X$ and consider  a germ of a smooth (everywhere below this means
$C^{\infty}$) ${\cal CR}$ mapping ${}_p{\bf f} : X \longrightarrow Y$ of $X$ to $Y$. 
 This means that there exists a neighborhood $U$ of $p$ in $\cc^n$ and a
representative mapping $f$ of ${}_p{\bf f}$ defined on $X \cap U$ such that $f(X
\cap U) \subset Y$.

  Consider also the field ${\cal M}_p(X)$  of restrictions to $X$
of germs of meromorphic functions  at $p$   and the finite type extension
${\cal M}_p(X)({}_p{\bf f}_1,...,{}_p{\bf f}_N)$  of this field generated by
the components of the germ  ${}_p{\bf f}$. Denote by $tr.deg. ({}_p {\bf f})$ the
transcendence degree  of this  extension over ${\cal M}_p(X)$ (see section 2
for definitions).  The transcendence degree  measures
"the degree of non - analyticity" of ${}_p{\bf f}$. If $tr.deg. {}_p{\bf f}$ is
equal to $m$, we will show that the graph of ${}_p{\bf f}$ is contained in a
complex $(n + m)$ - dimensional variety in $\cc^{n + N}$. In particular,   we will
show that ${}_p{\bf f}$ is real analytic if and only if $tr.deg.{}_p{\bf f} = 0$. 

Let $f: X \longrightarrow \cc^N$ be a smooth ${\cal CR}$ mapping of a real analytic
minimal hypersurface. There exists an open dense subset of $X$ where the rank of $f$
achieves its maximal value, which we call the generic rank of $f$. As usual, by
$rank {}_p{\bf f}$ (the rank  of the germ ${}_p{\bf f}$) we mean a generic rank of its
representative mapping; this definition does not depend on the choice of a
representative. The  geometric property of real submanifolds in $\cc^n$
which we study in this paper can be described as follows.

\begin{e-def}
\label{defA}
Let $Y$ be a real algebraic subset of $\cc^N$ and $p$ be a point in $Y$. Let also $r$
and $m$ be positive integers. We say that $Y$ is $(r,m)$ - flat at $p$ if there exists
a real analytic  submanifold $M \subset Y$  through $p$ of real dimension $r$ which
is  biholomorphic to the cartesian product $\Gamma \times D$, where $\Gamma$ is a real
analytic submanifold in a complex space of  smaller dimension  and $D$ is a domain
in $\cc^m$.  
\end{e-def} 

Of course, this definition still makes sense, if $Y$ is just real analytic. In this
paper we will use it in the algebraic case only.

One can consider the manifold $M$ in Definition \ref{defA} (up to biholomorphic
equivalence) as a trivial fiber bundle with {\it holomorphic} fibers and the real
analytic manifold $\Gamma$ as the base.  If $Y$ admits  such a bundle, holomorphically
embedded as a  submanifold, then in general one can  construct a smooth ${\cal CR}$
mapping from $X$ to $Y$ which is not real analytic. For instance, if a
manifold $\Gamma\times D$, where $\Gamma$ is a real submanifold in $\cc^k$ and $D$ is
a domain in $\cc$, is contained in $Y$ and $X$ is a real sphere, one can consider a
${\cal CR}$ mapping from $X$ to $\Gamma \times D$ of the form $(g,h)$, where $g: X
\longrightarrow \Gamma$ is a constant mapping and  $g$ is smooth, but not real
analytic ${\cal CR}$ function on $X$.

The main goal of the paper is to prove the following

\bigskip

{\bf Main Theorem} {\it Let $X \subset \cc^n$ ($n > 1)$ be a  real analytic
hypersurface, minimal at a point $p \in X$  and $Y\subset
\cc^N$ be  a real algebraic set. Let ${}_p{\bf f}: X \longrightarrow Y$
be a germ of a smooth ${\cal CR}$ mapping at $p$ with $rank {}_p{\bf f} = r$ ,
$tr.deg. {}_p{\bf f} = m$.  Then  in any neighborhood of ${}_p{\bf f}(p)$ there are
points where $Y$ is $(r,m)$ - flat}.

\bigskip

This result means that the existence of bundles described by Definition \ref{defA}
is the only reason   for the non - analyticity of a ${\cal CR}$ mapping and, moreover,
the "degree of non - analyticity" does not exceed the complex dimension of fibers.
 In order to use this theorem for the study of analytic properties of ${\cal CR}$ 
mappings it suffices to find the obstacles for a holomorphic embedding in $Y$ of 
such bundles. This gives an upper estimate for the transcendence degree $tr.deg.
({}_p{\bf f})$.  The following result is a corollary of Main Theorem (see section
4).

\begin{e-theo}
\label{theoB}
Let $X \subset \cc^n$ ($n > 1)$ be a real analytic hypersurface, $Y \subset \cc^N$ be
a real algebraic set  and  ${}_p{\bf f}: X \longrightarrow Y$ be a germ of a smooth
${\cal CR}$ mapping at $p \in X$. Suppose that $X$ is minimal at 
$p$ and denote by  $d$  the maximal dimension of complex analytic varieties
contained in $Y$ in a neighborhood of ${}_p{\bf f}(p)$. Then 
$tr.deg. ({}_p{\bf f})\leq d$ and the graph of ${}_p{\bf f}$ is contained in a
complex $(n + d)$ - dimensional analytic variety. In particular, if $Y$
does not contain complex analytic varieties of positive dimension in a
neighborhood of ${}_p{\bf f}(p)$, then ${}_p{\bf f}$ is real analytic.  
\end{e-theo} 

The estimate of the transcendence degree obtained in Theorem \ref{theoB} may be
considered as  a result on partial analyticity of ${}_p{\bf f}$. If, additionally, $X$
is algebraic, a  result similar to Theorem \ref{theoB} was obtained in \cite{CoMeSu};
however, there is an important difference. In \cite{CoMeSu} the transcendence degree
is considered over the field of {\it rational} functions in $\cc^n$. Hence, the
estimate on the transcendence degree concerns {\it algebraic} properties of $f$. The
proof of our results is more delicate since we assume $X$ to be only real analytic
and purely algebraic methods are not enough. Analytic part of this paper is based on
the ideas of \cite{Pi1}. 

Even in the case when $Y$ does not contain complex analytic varieties of positive
dimension, Theorem \ref{theoB} is new. We mention here some special cases of
Theorem \ref{theoB}  which were previously considered by other authors: 
\begin{itemize}
\item[1)] $X$ is a real analytic strictly pseudoconvex hypersurface in $\cc^n$ and
$Y$ is a real sphere in $\cc^N$ (\cite{Pi1}).
\item[2)] $X$ is a real analytic strictly pseudoconvex hypersurface in $\cc^n$ and $Y$
is a real algebraic set in $\cc^N$ with $d = 0$(A.Pushnikov \cite{Pu}). 
\item[3)]  $X$ and $Y$ are real algebraic strictly pseudoconvex hypersurfaces (
X.Huang \cite{Hu}).  
\item[4)]  $X$ is a real algebraic    hypersurface in $\cc^n$
containing no complex analytic varieties of positive dimension and $Y$ is the
sphere in $\cc^{n+1}$ ( M.S.Baouendi, X.Huang and L.P.Rothschild \cite{BaHuRo}).
\end{itemize}

We point out that the result of Theorem \ref{theoB} is  new even in the
equidimensional case. In particular, we have the following 

\begin{e-cor}
\label{corE}
Let $X$ and $Y$ be real hypersurfaces in $\cc^n$. Assume that $X$ is real analytic
and minimal and $Y$ is real algebraic. Let $ f: X \longrightarrow Y$ be  
a $C^{\infty}$ smooth ${\cal CR}$ mapping. If  $Y$ contains no complex
analytic varieties of positive dimension, then $f$ is real analytic.   
\end{e-cor}

In the  case when $X$ is real algebraic and holomorphically non - degenerate, a
similar result was obtained by M.Baouendi, X.Huang and L.Rothschild \cite{BaHuRo}.

Other  applications of Main Theorem  concern germs of ${\cal CR}$ mappings of
 maximal rank. This   assumption allows to  weaken the restrictions  on $Y$. 

  Recall that a real hypersurface $M$ is called holomorphically degenerate at a point
$a \in M$  if there exists a non-zero germ of a holomorphic vector
field tangent to $M$ at $a$; $M$ is called holomorphically non - degenerate in a
neighborhood $U$ of $p \in M$, if it is not holomorphically degenerate at every point
of $U \cap M$. Clearly,  $M \cap U$  is holomorphically non - degenerate, if and
only if in any  neighborhood of every point $ a \in M \cap U$ it is not locally
biholomorphic to the cartesian product $\Gamma \times D$ where $\Gamma$ is a real
analytic hypersurface a complex space of smaller dimension and $D$ is a domain in
$\cc$. One can
naturally extend this definition as follows. Let $M$ be a real analytic manifold in a
neighborhood of a point $p \in \cc^n$. We say that $M$ is $d$-holomorphically non -
degenerate in a neighborhood $U$ of $p$, if in any neighborhood of any point of $M
\cap U$ it is not biholomorphic to a cartesian product $\Gamma \times D$ where $D$
is a domain in $\cc^d$. Notice that a $d$ - holomorphically  non - degenerate
manifold can contain complex analytic varieties of dimension $\geq d$.

\begin{e-cor}
\label{corG}
Let $X \subset \cc^n$ be a real analytic hypersurface minimal at a point $p \in X$, 
 and $Y$ be a real algebraic  manifold in  $\cc^N$, $d$ -holomorphically non -
degenerate in a neighborhood of a point $p' \in Y$. Let ${}_p{\bf f}:
X\longrightarrow Y$   be a germ of a   smooth ${\cal CR}$ submersion with ${}_p{\bf
f}(p) = p'$. Then $tr.deg.({}_p{\bf f}) < d$ and the graph of ${}_p{\bf f}$ is
contained in a complex analytic variety of dimension $< n + d$.
\end{e-cor}    
 
In particular, we have

\begin{e-cor}
\label{corH}
Let $X$ and $Y$ be real hypersurfaces in $\cc^n$. Assume that $X$ is real analytic
and minimal at $p \in X$ and $Y$ is real algebraic and holomorphically non -
degenerate. Let $f: X\longrightarrow Y$ be a smooth ${\cal CR}$ mapping of 
generic rank $n$.  Then $f$ is real analytic.   
\end{e-cor}

If both $X$ and $Y$ are algebraic and holomorphically nondegenerate, a similar result
was obtained by M.S.Baouendi, X.Huang and L.P.Rothschild \cite{BaHuRo}. By the
example of P.Ebenfelt (see \cite{BaHuRo})   the minimality condition in  Corollary
\ref{corH} cannot be replaced by the holomorphic
nondegeneracy. Another special case of Corollary \ref{corH} ( $Y$ is a rigid
algebraic hypersurface) was considered by N.Mir \cite{Mir} and J.Merker and F.Meylan
\cite{MeMe}.  Related results were also obtained in \cite{MaMe,Me1}.

\section{Outline of the proof of the Main Theorem}

 Recall some standard  constructions from commutative algebra.  For any entire ring 
$K$ denote by  $K[x_1,\dots,x_m]$  the ring of polynomials in $m$ variables
with coefficients in $K$. If $F$ is the quotient field of $K$, we denote by  
$F(x_1,\dots,x_m)$ the quotient field of $K[x_1,...,x_n]$ which is identified with
the field  of rational fractions over $F$. In the case where $K = \cc$
we identify them with the ring $\cc[Z_1,\dots,Z_m]$ (resp. the field
$\cc(Z_1,\dots,Z_m)$) of polynomial (resp. rational) functions on $\cc^m$. More
generally, let $F \hookrightarrow E$ be  a field extension and
$\alpha_1,...,\alpha_n$ be elements of $E$. We denote by
$F(\alpha_1,...,\alpha_n)$ the smallest subfield of $E$, containing the field $F$ and
every $\alpha_j$. Obviously, it coincides with the field of rational fractions in
$\alpha_1,...,\alpha_n$ with coefficients in $F$. We say that $E$ is finitely
generated over $F$ if there exists a finite family of elements
$\alpha_1,...,\alpha_n$ in $E$ such that $E = F(\alpha_1,...,\alpha_n)$. Recall that
a finitely generated extension   
$E = F(\alpha_1,...,\alpha_n)$ is algebraic over $F$ if and only if every $\alpha_i$
is algebraic over $F$. If $S$ is a subset of $E$, we still denote by $F(S)$ the
smallest subfield of $E$ containing $S$ and say that $S$ generates $E$.

Let $F \hookrightarrow E$ be a field extension and $S$ be a subset of $E$. We say
that $S$ is algebraically independent over $F$ if for every $s_1,...,s_m
\in S$ and every polynomial $P \in F[x_1,...,x_m]$ the equality $P(s_1,...,s_m) =
0$ implies $P = 0$. A finite subset $S \in E$, which is algebraically independent over
$F$ and maximal with respect to the inclusion ordering, is called a transendence base
of $E$ over $F$. In this case $E$ is algebraic over $F(S)$. Any two transcendence
bases of $E$ have the same cardinality which is called the degree of transcendence of
$E$ over $F$ and is denoted by $tr.deg. (E/F)$. Moreover, if $S$ generates $E$ and
$S'$ is a subset of $S$, algebraically independent over $F$, then there exists a
transcendence base $B$ of $E$ such that $S' \subset B \subset S$. We refer the
reader to \cite{La} for  proofs of the above statements.

We describe now the general scheme of the proof of the Main Theorem.
Let $X$ be a generic real analytic manifold in a neighborhood of a point $p \in X$ in
$\cc^n$ and $Y$ be a real analytic subset in a neighborhood of a point $p' \in
\cc^N$.  We denote by ${\cal O}_p$ the ring of germs of holomorphic functions at $p$
and by ${\cal O}_p(X)$ the ring of germs of their restrictions to $X$ (since $X$ is
generic, these rings are isomorphic). Let ${\cal M}_p$ be the field of germs of
meromorphic functions at $p$ (i.e. the quotient field of ${\cal O}_p$) and ${\cal
M}_p(X)$ the field of germs of their restrictions on $X$. If $D$ is a domain in
$\cc^n$, ${\cal O}(D)$ denotes the ring of holomorphic functions on $D$ and ${\cal
M}(D)$ the field of meromorphic functions on $D$ i.e. the field of sections of the
sheaf of germs of meromorphic functions on $D$. Similarly,
${\cal O}(X)$ (resp. ${\cal M}(X)$) denote the ring (resp. field) of
their restrictions to $X$.  For every function $h$, meromorphic on $\Omega$, there
exists a complex analytic set in $\Omega$ such that $h$ is holomorphic outside of
this set. The smallest analytic subset of $\Omega$ with this property is called the
set of singular points of $h$ and is denoted by $Sing_h$. It is well known that
every meromorphic function on a Stein domain $\Omega$ is as a quotient of two
holomorphic functions in $\Omega$.

We denote by ${}_pC^{\infty}_{\cal CR}(X)$ the ring of germs of $C^{\infty}$ smooth
${\cal CR}$ functions on $X$ at $p$ and by $({}_pC^{\infty}_{\cal CR}(X))^N$ its
$N-th$ cartesian power. By a germ of a smooth ${\cal CR}$ mapping ${}_p{\bf f} : X
\longrightarrow Y$ we mean an element ${}_p{\bf f} \in   ({}_pC^{\infty}_{\cal
CR}(X))^N$ such that for some neighborhood $U$ of $p$ in $\cc^n$ there exists a
representative mapping $f$ of ${}_p{\bf f}$, smooth and ${\cal CR}$ on $X \cap U$
and satisfying $f(X \cap U) \subset Y$.  As in the previous section, we consider the field extension ${\cal
M}_p(X)({}_p{\bf f}_1,...,{}_p{\bf f}_N)$ of the field ${\cal M}_p(X)$ of germs of
meromorphic functions at $p$ finitely generated by components of ${}_p{\bf f}$  and
denote by $ tr.deg. ({}_p{\bf f})$ its transcendence degree over ${\cal M}_p(X)$.

  This means that for $m = tr.deg. ({}_p{\bf f})$ there exist integers $1\leq i_1
<\dots < i_m\leq N$ such that $({}_p{\bf f}_{i_1},\dots, {}_p{\bf f}_{i_m})$ is a
transcendence base of ${\cal M}_p(X)({}_p{\bf f}_1,\dots,{}_p{\bf f}_{N})$ over ${\cal
M}_p(X)$. After an eventual renumeration we  assume without loss of generality
that $i_1 = 1,\dots, i_m = m$. We will use the notation $ {}_p{\bf f} = ( {}_p{\bf g},
 {}_p{\bf h})$, where ${}_p{\bf g} = ({}_p{\bf g}_1,\dots, {}_p{\bf g}_m) =  ({}_p{\bf
f}_1,\dots, {}_p{\bf f}_m)$ and ${}_p{\bf  h} = ( {}_p{\bf h}_1,\dots, {}_p{\bf
h}_{N-m}) = ({}_p{\bf f}_{m+1},\dots, {}_p{\bf f}_{N})$.  Denote the coordinates in
$\cc^N$ by $Z' = (z',w')$, $z'= (z'_1,\dots,z'_m)$ and  $w' = (w'_1,\dots,w'_{N-m})$.

 Fix a representative mapping $f = (g,h)$ of ${}_p{\bf f}$ in a small enough connected
neighborhood $U$ of $p$ in $X$. The first step of our construction is the embedding
of the graph $\Gamma_{f}$ to  some complex variety, canonically defined by ${}_p{\bf
f}$. 

Since $\{{}_p{\bf g}_1,...,{}_p{\bf g}_m \}$ is an algebraically
independent system over the ring ${\cal M}_p(X)$, the substitution morphism  ${\cal
O}_p(X)[x_1,...,x_m]\longrightarrow {\cal O}_p(X)[{}_p{\bf g}_1,...,{}_p{\bf g}_m]$
is a ring isomorphism. By the definition of a transcendence basis and the Gauss
lemma there are non-zero irreducible polynomials $Q_j$ in  the ring ${\cal
O}_p(X)[{}_p{\bf g}_{1},\dots, {}_p{\bf g}_{m}][x]$ such that $Q_j({}_p{\bf h}_j) =
0$  (see \cite{La}, Ch.V, Th.10). Represent every polynomial $Q_j$  in the form
$Q_{j} =\sum_{k=0}^{N_j} {}_p{\bf q}_{jk}(Z, {}_p{\bf g})x^{k}$, where $deg Q_j = N_j
> 0$ and  ${}_p{\bf q}_{jk}
\in {\cal O}_p(X)[z']$. Denote by $q_{jk}$ corresponding representatives on
$U \times W$, where $W$ is a neighborhood of $g(p)$ in $\cc^m$. We associate with
them  holomorphic functions   

\begin{eqnarray}
\label{5.1'}
\hat{Q}_j(Z,Z') := \sum_{k = 0}^{N_j}q_{jk}(Z,z')w{'}_{j}^{k}
\end{eqnarray}
and   a complex  analytic variety ${\cal A}_{{}_p{\bf f}}$ in a
neigborhood $U \times U'$ of $(p,f(p))$ in 
$\cc^{n} \times \cc^{N}$, 
defined as the set of their common zeros:
\begin{eqnarray}
\label{5.1}
{\cal A}_{{}_p{\bf f}} = \{ (Z,Z'): \hat{Q}_j(Z,Z') = 0, j = 1,...,N-m \}
\end{eqnarray} 
 where every $\hat{Q_j}$ represents an irreducible polynomial
in $w'_j$ over the ring ${\cal O}_p(X)[z'_1,...,z'_m]$. 

 The graph $\Gamma_{f}$ of $f$ is contained in ${\cal
A}_{{}_p{\bf f}}$.  Let   
$$\lambda : \cc^{n}(Z) \times
\cc^{m}(z')\times\cc^{N-m}(w')\longrightarrow \cc^{n}(Z) \times \cc^{m}(z')$$ 
be the natural projection. Consider a  complex analytic variety ${\cal V}$ in a
neighborhood $U \times W$ of the point $(p,g(p))$ in $\cc^n \times \cc^{m}$ defined 
by $\cup_s \{ D_{s}(Z,z') = 0 \}$, where $D_s \in {\cal O}(U)[z']$ is the discriminant
of $\hat{Q}_s$ with respect to $w'_j$. Then the restriction
$\lambda : {\cal A}_{{}_p{\bf f}}\backslash
\lambda^{-1}({\cal V})\longrightarrow \cc^{n}\times \cc^{m}$ is a local
biholomorphism.  Denote by $\Gamma_{{}_p{\bf g}}$ the graph of 
$g$  in $U \times W$. Since every $D_s$ is an element of  ${\cal
O}(U)[z'_1,...,z'_m]$, the algebraic independence of the system $\{ {}_p{\bf
g}_1,\dots, {}_p{\bf g}_m\}$ over ${\cal O}_p(X)$ implies that $\Gamma_{{}_p{\bf g}}$
is not contained in ${\cal V}$.  Consider the set 

\begin{eqnarray}
\label{sing}
\Sigma = \cup_s \{ Z \in X\cap U: D_s(Z,g(Z)) = 0 \}
\end{eqnarray} 

Then for any $a \in X \backslash\Sigma$ the point
$(a,g(a))$ is not in ${\cal V}$. We will show in the next section (Lemma 3.1) that
$\Sigma$ is of measure 0 (with respect to $X$). 

\bigskip

{\bf Remark.} In general the variety ${\cal A}_{{}_p{\bf f}}$ is reducible. Let
$\tilde\Gamma_f = \{ (Z,Z') : Z \in X \backslash \Sigma, Z' = f(Z) \}$.  Since
$\lambda : {\cal A}_{{}_p{\bf f}} \longrightarrow U \times V$ is a local
biholomorphism near any $(Z,Z') \in \tilde\Gamma_f$, the union of all irreducible
components of ${\cal A}_{{}_p{\bf f}}$  with non - empty intersections with
$\tilde\Gamma_f$ is a complex analytic set in $U \times U'$ of pure dimension  $n +
m$,  containing the graph  $\Gamma_f$ (one can easily show that, in fact, there is
only one such component, but we do not need this fact).

\bigskip

The next step is  to study  the local geometry of
${\cal A}_{{}_p{\bf f}}$ in a neighborhood of  a point $a
\in (X\cap U) \backslash \Sigma$. Let $\pi: \cc^n \times \cc^N \longrightarrow \cc^n$
and $\pi':\cc^n\times\cc^N\longrightarrow \cc^N$ be the natural projections. We denote
by ${\cal A}_{{}_p{\bf f}} \vert X$ the "restriction" of ${\cal A}_{{}_p{\bf f}}$ to
$X$ defined by the additional assumption $Z \in X$ in equations (\ref{5.1}) i.e.
${\cal A}_{{}_p{\bf f}} \cap \pi^{-1}(X)$. Obviously, $\pi({\cal A}_{{}_p{\bf f}}\vert
X)\subset X$.

It follows from the implicit function theorem that there exist a neighborhood $U_a$
of $a$ in $\cc^n$, a neighborhood $V'_a$ of $g(a)$ in $\cc^m$ and  a neighborhood
$V''_a$ of $h(a)$ in $\cc^{N-m}$  such that ${\cal
A}_{{}_p{\bf f}}
\vert X$ can be represented in $U_a \times V_a$ with $ V_a = V'_a \times V''_a$ in
the form
$w' = H(Z,z')$, $Z \in X \cap U_a$, where $H$ is holomorphic in $U_a \times V'_a$.
Hence the triple $\{ {\cal A}_{{}_p{\bf f}}\vert X \cap (U_a
\times V_a), \pi, X \cap U_a \}$ has the structure of a trivial fiber bundle over
$X \cap U_a$ with holomorphic fibers $\pi^{-1}(Z) = \{ Z' = (z',w'): w' = H(Z,z')
\}$ of complex dimension $m$.

The crucial question  is whether the projection $\pi'$ takes the fibers of
$\pi$ to $Y$? Thus  we need  the following condition :

\begin{eqnarray}
\label{supsup5.1}
\pi'(({\cal A}_{{}_p{\bf f}}\vert X) \cap (U_a \times V_a) \subset Y  
\end{eqnarray}

This condition  allows to transfer
the structure of the bundle $\{ ({\cal A}_{{}_p{\bf f}} \vert X) \cap (U_a \times
V_a),\pi, X\cap U_a\}$ to the target set $Y$. Since the fibers of $\pi$ are complex
manifolds and the projection $\pi'$ is holomorphic in the ambient space, it preserves
the complex structure of fibers. Moreover, the restriction of $\pi'$ to a fiber of
$\pi$ is biholomorphic and this allows to construct many complex subvarieties in $Y$.

The following proposition shows that the condition (\ref{supsup5.1}) implies the
$(r,m)$ - flateness of $Y$ and is the key for our approach.

\begin{e-pro}
\label{prosup5.1}
Suppose that (\ref{supsup5.1}) holds. Then  there exist a positive integer $r = r(a)
\geq rank {}_p{\bf f}$ such that in any neighborhood of $f(a)$ there are points where
$Y$ is $(r,m)$ - flat.
\end{e-pro}
\proof This is a direct corollary of the rank theorem. Denote by $r$ the maximal rank
of the restriction $\pi'\vert ({\cal A}_{{}_p{\bf f}} \vert X) \cap
(U_a \times V_a)$ and by $S$ the open dense subset of $({\cal A}_{{}_p{\bf f}}
\vert X) \cap (U_a \times V_a)$ where this rank is equal to $r$.
Since the graph of $f$ over $X$ is contained in ${\cal A}_{{}_p{\bf f}}\vert X$,
$r$ is bigger or equal to  the generic rank of $f$.  Fix a point $(Z^0,Z'^0) \in
S$, close enough to $(a,f(a))$, and neighborhoods $W_a \subset U_a$ of
$Z^0$ in $\cc^n$ and $W'_a\subset V_a$ of $Z'^0$ in $\cc^{N}$, such that the
intersection $\Omega = ({\cal A}_f \vert X)\cap (W_a\times W'_a)$ is contained in
$S$. It follows from  (\ref{supsup5.1}) and the real analytic
version of the rank theorem, applied to the restriction $\pi'\vert \Omega$,  that
there exist  neighborhoods $\hat{W}_a \subset W_a$ of $Z^{0}$ in $\cc^n$ and
$\hat{W}'_a \subset W'_a$ of $Z'^0$ in $\cc^{N}$  such that $M_a := \pi'(\Omega)$ is
an $r$-dimensional real analytic submanifold in $Y \cap W'_a$. Furthemore, since
the restriction of the projection $\pi'$  on every fiber of $\pi$ has the
maximal rank, there exists a real analytic submanifold $\Gamma_a
\subset X$ through $Z^0$  such that the restriction $\pi':
\Omega\cap\pi^{-1}(\Gamma_a)\longrightarrow M_a$ is a real analytic diffeomorphism.
Since the projection $\pi'$ is holomorphic on the ambient space and 
$\Omega\cap\pi^{-1}(\Gamma)$ is a trivial fiber bundle over
$\Gamma_a$ with holomorphic fibers, we get desired statement. 

\bigskip

Since $(X \cap U) \backslash \Sigma$ is dense in $X \cap U$, this proves the Main
Theorem under the additional condition (\ref{supsup5.1}).
 
The goal of the rest of this paper is to show that the condition
(\ref{supsup5.1}) holds if $X$ is a real analytic minimal hypersurface in $\cc^n$ and
$Y$ is a real algebraic set in $\cc^N$. We make use of two principal tools: 
results on meromorphic extension of certain classes of ${\cal CR}$ functions and 
 algebraic field extensions. In the next section we discuss the problem of
meromorphic extension.

\section{Reflection, meromorphic extension and algebraic dependence over
intermediate fields}

{\bf 3.1. Intermediate field extensions.} Let $X$ be a  real analytic generic
minimal (in the sense of A.Tumanov \cite{Tu}) submanifold in a neighborhood of a point
$p \in X$ in $\cc^n$. Assume that $X$ is defined by a vector - valued function $\rho =
(\rho_1,...,\rho_d)$ which is real analytic in a neighborhood of $p$ and $\partial
\rho_1
\wedge ...\wedge \partial \rho_d \neq 0$ at $p$. Denote by $Z = (z,w)$ with 
$z = (Z_1,...,Z_{n-d})$  and $w = (Z_{n-d+1},...,Z_n)$ the coordinates in $\cc^n$.
Assume that  $det ( \frac{\partial \rho}{\partial w})(p) \neq 0$.  Consider the Cauchy
- Riemann operators in a neighborhood of $p$ on $X$ given by 

\begin{eqnarray}
\label{4.1}
{\cal L}_k = \frac{\partial}{\partial \overline
z_k} -  \sum_{j =
1}^{d}a_{jk}(Z)\frac{\partial}{\partial\overline w_j} 
\end{eqnarray}
with $(a_{jk}) = (\partial \rho/\partial w)^{-1}(\partial \rho/\partial w)$.

Denote by ${\cal O}_p^{\R}(X)$ the ring of germs at $p$ of real analytic functions on
$X$ and by   ${}_pC^{\infty}_{{\cal CR}}(X)$ (resp.  
${}_pC^{\infty}_{\overline{\cal CR}}(X)$) the ring of germs at $p$ of 
$C^{\infty}$ smooth  ${\cal CR}$ (resp. anti ${\cal CR}$) functions on $X$. If
$\Omega$ is an open subset of $X$, we denote by  ${\cal O}^{\R}(\Omega)$ (resp.
$C^{\infty}_{\overline{\cal CR}}(\Omega)$) the ring of real analytic (resp. anti
${\cal CR}$) functions on $\Omega$.

 We denote by ${\cal P}_p(X)$ the set of linear combinations of elements of 
 ${}_pC^{\infty}_{\overline{\cal CR}}(X)$ with coefficients in  ${\cal O}_p^{\R}(X)$.
 It is worthwhile to emphasize that the values of real analytic coefficients and anti
${\cal CR}$ functions  are always considered at the same point of $X$, so the elements
of ${\cal P}_p(X)$ are germs of functions on $X$ at $p$.   

Notice that ${\cal P}_p(X)$ is a ${\cal O}_p^{\R}(X)$ - module
and  ${\cal L}_j({}_p{\bf h}) \in {\cal P}_p(X)$ for every germ
${}_p{\bf h} \in {\cal P}_p(X)$.

We will need the following  property of the ring ${\cal P}_p(X)$.

\begin{e-lemme}
\label{lem4.1}
Let $X$ be a generic real analytic  manifold  in $\cc^n$ minimal at a point $p \in
X$, and ${}_p{\bf h} \in  {\cal P}_p(X)$. Let $U$ be  a connected
neighborhood  of $p$ in $\cc^n$,  $h$ be a representative   of ${}_p{\bf
h}$ on $X\cap U$ and $E = \{Z \in X \cap U:  h = 0 \}$. Suppose that   that the
intersection of
$E$ with every neighborhood of $p$ is of positive measure (with respect to $X$). 
Then ${}_p{\bf h} = 0$.
\end{e-lemme}

\proof  It follows from \cite{Tu} that there exist small enough connected
neighborhoods $U'
\subset U$ of $p$  and  an open convex cone ${\cal V}$ in $\R^d$ such that
each ${\cal CR}$ function on $X \cap U$ extends holomorphically to the wedge ${\cal
W} = \{ Z \in U' \vert r(Z) \in {\cal V}\}$.    

Consider a foliation of $X \cap U'$ by real analytic maximal totally real manifolds
$\{ M_t \}_t$ (i.e. $dim M_t = n$ for any $t$), where the parameter $t$ is in a
domain $D \subset \R^{n-d}$.    By the Fubini - Tonelli theorem  $F:= \{ t \in D
: mes (M_t \cap E) > 0 \}$ is the subset of $D$ of positive measure. 
Let $t \in F$ and $h = \sum_{\nu = 1}^{s} a_{\nu} g_{\nu}$ with $a_{\nu} \in {\cal
O}^{\R}(X \cap U)$ , $g_{\nu} \in C^{\infty}_{\overline{\cal CR}}(X \cap U)$.

 Every restriction $a_{\nu} \vert M_{t}$   extends 
antiholomorphically  to a neighborhood of $M_t$ and  every $g_{\nu}$ extends
antiholomorphically to ${\cal W}$. Hence, $h \vert M_t$ coincides with some
function $\hat{h}$ which is antiholomorphic  ${\cal W}$ and smooth up to $M_t
\subset X \cap U$. By the boundary uniqueness theorem \cite{Co,Pi2}   
$\hat{h} \equiv 0$  and thus $h \vert M_t = 0$. Therefore, $\cup_{t \in F} M_t
\subset E$.  

Now we can consider another foliation $\{ M'_t \}$ of $X \cap U$ such that each $M'_t
\cap E$ has positive measure and repeat the previous arguments to conclude that
$E = X \cap U$. 

\bigskip

In particular, we have the following 

\begin{e-cor}
If $X$ is minimal at $p$, then ${\cal P}_{p}(X)$ is an entire ring.
\end{e-cor}

This last corollary allows us to introduce the quotient field of the ring ${\cal
P}_p(X)$ which we denote by ${\cal M}^{*}_p(X)$. Let $U$ be a neighborhood of $p$
in $\cc^n$. A representative $h$ of
${}_p{\bf h}$ of ${\cal M}^{*}_p(X)$ in $X \cap U$ is a quotient ${\varphi}/{\psi}$,
where $\varphi$ and $\psi$ are representatives of  germs from ${\cal P}_p(X)$ in $X
\cap U$.   Let  $Sing_{h}:=\{ Z \in X \cap U:\psi = 0 \}$ be the singular set of $h$. 
If $U$ is small enough, it is a closed subset  of $X \cap U$ of measure $0$ and   $h$
is a smooth function on $(X \cap U) \backslash Sing_h$.  We will call the points of
$(X
\cap U)\backslash Sing_h$ regular for $h$.  In what follows we say that a subset
$E \subset (X\cap U)$ through $p$ is a ${\cal P}$ - set, if there exists a non-
zero germ ${}_p{\bf g} \in {\cal P}_p(X)$  which admits a representative $g$ on $X
\cap U$ such that $E = \{Z \in X \cap U: g = 0 \}$.

\bigskip

{\bf 3.2. Reflection and meromorphic extension.} The following assertion is our main
analytic tool.

\begin{e-pro}
\label{pro2.2}
Let $X$ be a  real analytic hypersurface in $\cc^n$, minimal at $p \in X$.
 Let ${}_p{\bf h} \in {\cal M}^{*}(X)$ and $h$ be its representative in a
neighborhood $U$ of $p$. Assume that there exists a  ${\cal P}$ - set $E \subset (X
\cap U)$ containing $Sing_{h}$ such that $h$ is a ${\cal CR}$ function on $(X \cap U)
\backslash E$. Then ${}_p{\bf h} \in {\cal M}_p(X)$, i.e. ${}_p{\bf h}$ is
meromorphic.
\end{e-pro}

\proof In what follows we denote by $\Delta^k(a,r)$ the  polydisc in $\cc^k$
centerd at $a$ of radius $r$ (we write $\Delta(a,r)$ if $k = 1$ and $\Delta^k(r)$
if $a = 0$).  We may assume that $p = 0$, $U$ is a neighborhood of the origin in
$\cc^n$ of the form $U = \Delta^n(r)$, $r > 0$. In what follows we will replace $r$
by smaller positive numbers when we need; in order to avoid complications of the
notations, we keep the same notation $r$ for these numbers. Suppose also that  $X
\cap U$ is defined by  $\{ Z \in U: \rho(Z) = 0 \}$ where $\rho$ is a real
analytic function with $\partial \rho/\partial Z_n \neq 0$ on $U$ and $X$ is
minimal at every point of $X \cap U$. 

In view of Lemma \ref{lem4.1} we may assume that $E \cap U$ is of
measure $0$ (with respect to $X$). We  use the notation $Z = (z,w)$, $z \in
\cc^{n-1}$, $w \in \cc$ for coordinates in $\cc^n$ and  assume that for every $z
\in \Delta^{n-1}(r)$ the linear disc 

\begin{eqnarray}
\label{3.3}
l(z): = \{ z \} \times \Delta(r)
\end{eqnarray} 
intersects $X$ transversally at every point. 
  
If $V \subset U$ is a neighborhood in $\cc^n$ of a point $a \in (X \cap U)$, $V^+$
(resp. $V^-$) denotes the corresponding one - sided neighborhood $\{ Z \in V: \rho(Z)
> 0\}$ (resp. $\{ Z \in V: \rho(Z) < 0 \}$).

 In view of \cite{Tr} we may suppose that the following holds:  

\begin{itemize}
\item[(a)]  there exist a fundamental system of neighborhoods $\{ U_s \}_s$, $\{ V_s
\}_s$  of the origin  with $ V_s \subset U_s$ and such that for any 
$s$ every function holomorphic in $U_s^{+}$ extends holomorphically to $V_s$.
\item[(b)] for every point $a \in X \cap U$   there exist fundamental systems 
of neighborhoods $\{{}_aV_s \}_s$, $\{ {}_aU_s \}_s$ , ${}_aV_s \subset {}_aU_s$ of
$a$ such that  for any $s$ every function holomorphic in ${}_aU_s^{+}$ extends
holomorphically to ${}_aV_s$ or  for any $s$ every function holomorphic in
${}_aU_s^{-}$ extends holomorphically to ${}_aV_s$. 
\end{itemize}

  Let $h = \varphi/\psi$, where the functions $\varphi$ and $\psi$ are of  the form 

\begin{eqnarray}
\label{*}
\sum_{\nu = 1}^{s} a_{\nu} g_{\nu} 
\end{eqnarray}
with $a_{\nu} \in {\cal O}^{\R}(X \cap U)$ , $g_{\nu} \in C^{\infty}_{\overline{\cal
CR}}(X \cap U)$. We proceed the proof in three steps.

\bigskip

{\it Step 1.  Reflection.} We begin with the following corollary of the
Schwarz reflection principle. 

\begin{e-lemme}
\label{lem2.3}
Let $\phi^{-}$ be a function, holomorphic on  $U^{-}$ and continuous up to 
$X \cap U$. Then there exists  a real analytic function $\phi^{+}$ on $U^{+}$
continuous up to  $X \cap U$ and such that for every $z \in \Delta_{n-1}(r)$ the 
restriction $\phi^{+}\vert l(z) \cap U^+$  is antiholomorphic
and $\phi^{-} \vert (X \cap U)  = \phi^{+} \vert X (\cap U)$. 
\end{e-lemme}

\proof It follows from the implicit function theorem that there exists
 a real analytic mapping $\Phi(z,w): U \longrightarrow \Delta$ with the
following property:  for every $z \in \Delta^{n-1}(r)$ there exists $\delta_z >
0$ such that the mapping $\Phi_z: \Delta(r) \longrightarrow \Delta(\delta_z)$, $\Phi_z
: w \mapsto \Phi(z,w)$ is a biholomorphism between the domains  $\{ w \in \Delta(r):
r(z,w) < 0 \}$ and $\{ \omega \in \Delta(\delta_z): Im \omega < 0
\}$ and maps the real analytic curve
$\{ w\in\Delta: r(z,w) = 0 \}$ to $[-\delta_z,\delta_z]$. Let $\sigma: \omega
\longrightarrow \overline\omega$ be the standard conjugation in $\cc$. Now it suffices
to set $\phi^{+}(z,w) = \phi^{-} \circ \Phi_z^{-1} \circ\sigma \circ
\Phi_z(w)$.  

\bigskip

Every real analytic coefficient in the representation (\ref{*}) of $\varphi$ (resp.
$\psi$) extends  in a real analytic way to $U^+$ and is holomorphic on every disc
(\ref{3.3}). By the condition (a) every anti ${\cal CR}$ function on $X \cap U^+$
extends antiholomorphically to a smaller neighborhood of the origin which we again
denote by $U$. Apply Lemma \ref{lem2.3} to these antiholomorphic functions in the
representation (\ref{*}) of $\varphi$ and $\psi$. We get  that they extend 
to $U^+$ as real analytic functions,  holomorphic on every disc $l(z)
\cap U^+$. Finally,  we get that $\varphi$ and $\psi$ extend to  real analytic
functions $\varphi^{+}$ and $\psi^{+}$ on $U^{+}$, which are holomorphic on every
disc $l(z) \cap U^{+}$.

   Let $E'$ denote the set of points $z$ in $\Delta^{n-1}(r)$ such that
the disc $l(z)$ intersects $E$ along a set of positive linear measure. Since
every function of the form (\ref{*}) extends antiholomorphically to $l(z) \cap
U^-$, it follows from  the boundary uniqueness theorem that for every $z
\in E'$,  one has $(l(z)\cap X)\subset E$. By  Fubini - Tonelli theorem, $E'$
is a closed subset of $\Delta^{n-1}(r)$ of measure $0$. Hence,  $h^{+} :=
\varphi^{+}/\psi^{+}$ is the quotient of two real analytic in $U^{+}$ functions  and
for every $z \in \Delta^{n-1}(r) \backslash E'$ the restriction $h^{+}
\vert l(z) \cap U^{+}$ is  meromorphic  and coincides with $h$ on $l(z) \cap
((X \cap U) \backslash E)$. 

\bigskip

{\it Step 2.  Meromorphic extension.} The
next step is to extend $h$ meromorphically to an open dense subset of $U^+$. 
Let $S = \{ Z = (z,w) \in U^{+}: Z \in l(z), z \in E' \}$.

\begin{e-lemme}
\label{suplem2.4}
The function $h^+$ is meromorphic on $U^{+} \backslash S$.
\end{e-lemme}

\proof Let $(z_0,w_0)$ be a point in $U^{+} \backslash S$.  Fix
a point $a \in l(z_0) \cap ((X \cap U)
\backslash E))$. By  (b) there exists a one - sided neighborhood $V$ of
$a$ such that $h$ extends holomorphically to $V$. Since $h^+$ is meromorphic on every
$l(z)$ and coinsides with the ${\cal CR}$ function $h$ on $(X\cap U) \backslash E$, it
 is holomorphic on $V$, if $V \subset U^+$, or extends holomorphically to $V$, if $V
\subset U^-$. Let us assume that $V \subset U^+$ (the other case can be treated
similarly).

Fix a point $\omega \in \Delta(r)$ with $(z_0,\omega) \in V$ and
$\delta' > 0$ such that $\{z_0 \} \times  \Delta(\omega,\delta')$ is contained in
$V$. There exists a simply connected domain $G$ in $\cc$ such that $\{ z_0 \}\times
G$ is compactly contained in $U^+ \backslash E$ and contains both $\{ z_0 \} \times
\Delta(\omega,\delta')$ and $(z_0,w_0)$. Fix $\delta'' > 0$ such that the polydisc
$\Delta^{n-1}(z_0,\delta'') \times \Delta(\omega,\delta')$ is contained in $V$. The
function $h^+$ is holomorphic there and for every fixed $z
\in \Delta^{n-1}(z_0,\delta'')$ is meromorphic on $\{ z
\} \times G$. It follows now from the classical Rothstein lemma \cite{Ro,Sh}
that $h^+$ is meromorphic on $\Delta^{n-1}(z_0,\delta'')
\times G$. Thus, $h^+$ is meromorhic in a neighborhood of any point of $U^{+}
\backslash S$.

\bigskip

{\it Step 3. Elimination of singularities.}   We show now that $S$ is a removal
singularity for $h^+$.

\begin{e-lemme}
\label{lemsupsup3.6}
Let $D \subset \cc$ be a domain and $\varphi$, $\psi$ be real analytic functions in
$D$. Suppose that $\psi(\zeta^0) \neq 0$ for some $\zeta^0 \in D$ and $h: =
\varphi/\psi$ is holomorphic in a neighborhood of $\zeta^0$. Then $h$ is meromorphic
in $D$.
\end{e-lemme}
\proof  It is enough to prove lemma for $D = \Delta = \{ \vert \zeta \vert < 1 \}$
and $\zeta^0 = 0$. Let $\zeta' \in \Delta \backslash \{ 0 \}$ and $\gamma = \{ \zeta
\in \Delta :\zeta = t \zeta', t \in [0, 1/\vert \zeta' \vert [ \}$. The restrictions
$\varphi\vert \gamma$ and $\psi \vert \gamma$ admit holomorphic extensions
$\hat\varphi$ and $\hat\psi$ to some neighborhood $\Omega$ of $\gamma$ in $\Delta$.
Therefore the function $\hat\varphi/\hat\psi$ is meromorphic in $\Omega$ and
coincides with $h = \varphi/\psi$ on $\gamma$ near $0$. By the uniqueness theorem $h
= \hat{\varphi}/\hat{\psi}$ in a neighborhood of $0$. By analyticity the equality
$\varphi\hat{\psi} = \hat{\varphi}\psi$ holds everywhere in $\Omega$ and therefore
$h$ is meromorphic in $\zeta' \in \Omega$.

\begin{e-lemme}
The function $h^{+}$ extends meromorphically to $U^{+}$.
\end{e-lemme}
\proof  Let $Z^0 \in S$ and $\tilde{E} = \{ Z \in U^+: \psi^+(Z) = 0 \}$. We may
assume that $Z^0 = 0$. Moreover, since $\tilde{E}$ is closed and nowhere dense in
$U^+$, we can choose an affine coordinate system $Z = (Z_1,'Z)$, $'Z =
(Z_2,...,Z_n)$ in $\cc^n$ with the following properties:
 
there exist a point $a = (a_1,'0) \in U^+ \backslash \tilde{E}$ and positive
$\delta_1 < \delta_2$ such that 
\begin{itemize}
\item[(i)] $U_1 : = \Delta^n(a,\delta_1) = \Delta(a_1,\delta_1) \times
\Delta^{n-1}('0,\delta_1) \subset U^+ \backslash \tilde{E}$,
\item[(ii)] $U_2 : = \Delta(a_1,\delta_2) \times \Delta^{n-1}('0,\delta_1) \subset
U^+$,
\item[(iii)] $0 \in U_2$.
\end{itemize}

The function $h^+ = \varphi^+/\psi^+$ is already known to be holomorphic in $U_1$.
By Lemma \ref{lemsupsup3.6} for any fixed $'Z \in \Delta^{n-1}('0,\delta_1)$ it
extends meromorphically to $\Delta(a_1,\delta_2)$ as a function of $Z_1$. Hence, by
the Rothstein lemma, it is meromorphic in a neighborhood of $Z^0 = 0$.

\bigskip

It follows from (b) that the envelope of holomorphy of $U^{+}$ contains a neighborhood
of the origin. Therefore,$h^{+}$ extends meromorphically to a neighborhood of the
origin (see for instance \cite{Iv,Sh,Sh1}). This completes the proof of Proposition
\ref{pro2.2}.

\bigskip

{\bf 3.3. Algebraic dependence.} We  will now apply the previous results to
analyze the property of algebraic dependence over certain functional fields.

\begin{e-def}
Let $\{ {}_p{\bf h}_1, \dots, {}_p{\bf h}_k \}$ be a finite subset of
${}_pC^{\infty}_{\cal CR}(X)$. We say that it is  algebraically dependent over ${\cal
M}_p^{*}(X)$ if there exist a non-zero polynomial ${}_p{\bf Q}$ in ${\cal
M}_p^{*}(X)[X_1,\dots,X_k]$, a neighborhood $U$ of $p$, representative functions $h_j$
of ${}_p{\bf h}_j$ on $X \cap U$ , a representative polynomial $Q$ of ${}_p{\bf
Q}$ on $X \cap U$ and a ${\cal P}$ - set $E \subset X
\cap U$  containing the singular set  of every coefficient of $Q$, 
 such that $Q( h_1,\dots, h_k) = 0$ on $(X \cap U)\backslash E$.
\end{e-def}

The main result of this section  is the following:

\begin{e-pro}
\label{pro4.3}
If a  finite subset of ${}_pC^{\infty}_{\cal CR}(X)$ is algebraically dependent  over
${\cal M}_p^{*}(X)$, then it  is also  algebraically dependent over ${\cal M}_p(X)$.
\end{e-pro}

\proof Let $\{ {}_p{\bf h}_1, \dots,  {}_p{\bf h}_k \} \subset {}_pC^{\infty}_{\cal
CR}(X)$ be a algebraically dependent system  over ${\cal M}_p^{*}(X)$. This means that
there exist a neighborhood $U$ of $p$ , a system of representatives $\{ h_1,...,h_k
\}$ on $X \cap U$,  a  polynomial 
$$Q(X_1,...,X_k) = \sum_{j =0}^{m}q_{j}(Z,\overline{Z})S_{j}(X_1,\dots,X_k),$$
where $S_{j}$ are monomials and every $q_j$ represents a non - zero germ in ${\cal
M}_p^{*}(X)$, such that $Q( h_1,\dots, h_k) = 0$ on $(X\cap U)\backslash E$ . Here
$E \subset X \cap U$ is a ${\cal P}$ - set   containing the singular set of
every $q_j$. After   dividing $Q$ by $q_m$ we additionly have $q_{m}(Z,\overline{Z})
= 1$ and the degree of the monomial $S_m$ is not zero.  We will prove the statement by
induction in $m$.

For $m = 1$ we have $S_{1}( h_1,\dots, h_k) = 0$, so at least one of ${}_p{\bf h}_j$
is zero and the set  $\{ {}_p{\bf h}_1, \dots,  {}_p{\bf h}_k \}$ is
algebraically dependent over  ${\cal M}_p(X)$. 

Now  assume that the desired assertion is true for any polynomial containing
$\leq m - 1$ terms and apply  the Cauchy - Riemann operators  (\ref{4.1})  to
the equation  $Q(h_1,\dots, h_k) = 0$.
There are two possibilities
\begin{itemize}
\item[(a)] If ${\cal L}_{s}(q_{j}) = 0$ on $(X \cap U) \backslash
E$ for any $s$ and $j$. Then    Proposition \ref{pro2.2} implies  that the
coefficients
$q_j$  are meromorphic and  the set $\{ {}_p{\bf h}_1, \dots,  {}_p{\bf h}_k \}$ is
algebraically dependent over  ${\cal M}_p(X)$.  
\item[(b)] If there exists $j_{0}$ and $s_{0}$ such that ${\cal L}_{s_{0}}(q_{j_0})
\neq 0$ ,  we may apply the induction assumption to the polynomial $ \sum_{j =
0}^{m-1}({\cal L}_{s_0}q_{j})S_j$.
\end{itemize}   

\bigskip

We have the following 

\begin{e-cor}
\label{cor4.5}
Let ${}_p{\bf R} \in {\cal O}_p^{\R}(X)[w,\overline{w}]$ be a
polynomial in $w,\overline{w} \in \cc^m$ with real analytic (in $Z \in X$) 
coefficients and ${}_p{\bf g}= ( {}_p{\bf g}_1,...,{}_p{\bf g}_m)$ be a germ of a
mapping with components in ${}_pC^{\infty}_{\cal CR}(X)$ which are algebraically
independent over ${\cal M}_p(X)$. Suppose that there exist a neighborhood $U$ of $p$,
a representative mapping $g$ of ${}_p{\bf g}$ defined on $X\cap U$, a representaive
polynomial $R \in {\cal O}^{\R}(X \cap U)[w,\overline{w}]$ of ${}_p{\bf R}$ and a
${\cal P}$ -  set $E$ in $X \cap U$ such that $R(Z,g,\overline{g})$ vanishes on $(X
\cap U) \backslash E$. Then ${}_p{\bf R} = 0$ in ${\cal O}_p^{\R}(X)[w,\overline{w}]$.
\end{e-cor}

\proof  We represent $R$ in the form $R = \sum_{J}b_{J}(Z,\overline{w})w^{J}$ with
$b_{J} \in {\cal O}^{\R}(X \cap U)[\overline{w}]$ . Assume that  there exists a
coefficient $b_{J_{0}}$ such that $b_{J_{0}}(Z,\overline{ g})$ represents a non
-zero element in ${\cal M}_p^{*}(X)$. Hence, the system $\{ {}_p{\bf g}_1,\dots,
{}_p{\bf g}_m\}$ is algebraically dependent over ${\cal M}_p^{*}(X)$. Then Proposition
\ref{pro4.3} implies that this system is algebraically dependent over ${\cal
M}_p(X)$: this is a contradiction. Therefore,  for any $J$ the coefficient 
$b_{J}(Z,\overline{ g})$ represents zero in ${\cal M}_p^*(X)$. After the complex
conjugation we can write  $\overline{b_{J}} = \sum_{I}
c_{JI}(Z)w^{I}$. If there exist $ J_0 , I_0$ such that $c_{J_{0}I_{0}}$ represents
a non - zero element in ${\cal O}_p^{\R}(X)$, we apply  Proposition \ref{pro4.3}
again and obtain a contradiction as above. Thus every $c_{JI} = 0$ in ${\cal
O}^{\R}(X \cap U)$ and  $R(Z,w,\overline{w}) = 0$ in ${\cal O}^{\R}(X \cap
U)[w,\overline{w}]$. This completes the proof.

\section{ Completion of the proofs }

 We suppose  that we are in the settings of the Main Theorem  and will use the
notation of section 2.  It suffices  to show that (\ref{supsup5.1}) holds. 

Let the target real algebraic set $Y$ be defined by  real polynomials
$P'_{k}(Z',\overline{Z'})$ , $k = 1,\dots,q$, $Z' \in \cc^{N}$. 
Since $f$ takes $X$ to $Y$, we have $P'_{k}( f,\overline{ f}) = 0, k = 1,\dots,q$ .

In this section we use polarization techniques, so we consider real  analytic
functions as  functions in $Z,\overline{Z}$. Recall also that $U$ is a neighborhood
of $p$ in $\cc^n$, chosen in section 2 and $\Sigma$ is defined by (\ref{sing}).

\begin{e-pro}
\label{pro5.3}
There exist  polynomials $R_k \in {\cal O}^{\R}(X \cap U)[z',\overline{z}']$, $k =
1,...q$ with the following property: for every point $a  \in (X \cap U) \backslash
\Sigma$  there exist  neighborhoods $V_1$ of $(a,g(a))$ in $\cc^{n+m}$ and $V_2$ of
$h(a)$  in $\cc^{N-m}$ such that  the intersection $({\cal A}_{{}_p{\bf f}} \cap (X
\times Y)\cap (V_1 \times V_2))$  is defined by 
$$\{ (Z,z',w') \in V_1 \times V_2: Z\in X, \hat{Q}_j(Z,z',w_j) = 0,
R_k(Z,\overline{Z},z',\overline{z}') = 0, j = 1,...,N, k = 1,...,q\}$$
\end{e-pro}

\proof  The intersection ${\cal A}_{{}_p{\bf f}}\cap (X \times Y)$ is
defined by the equations 
$$\hat{Q}_j(Z,z',w'_j) = 0, P'_k(Z',\overline{Z}') = 0$$ 
$ Z\in X \cap U$. We can add the conjugate equations $\overline{\hat{Q}_j} = 0$ that
does not change the solutions.  Consider the complexified system 

\begin{eqnarray}
\label{sup5.8}
& &\hat{ Q}_j(Z,z',w'_j) = 0, {\tilde{Q}}_j(\zeta,\tau,\omega_j) = 0,
P'_k(z',\tau,w',\omega) = 0
\end{eqnarray} 
where $j= 1,...,N$, $k = 1,...,q$ , $Z = \overline{\zeta}$, $z' =\overline{\tau}$, 
$w' = \overline{\omega}$ and  ${\tilde{Q}}_j(\zeta,\tau,\omega_j):=
\overline{\hat{Q}}_j(\overline\zeta,\overline\tau,\overline\omega_j)$.

By the elimination theory \cite{W},\cite{Mu}  there exist a
neighborhood ${\cal U}$ of the point $(p,\overline{p})$ in $\cc^n(Z) \times
\cc^n(\zeta)$ and   functions $R_k(Z,z',\zeta,\tau) \in {\cal O}({\cal
U})[z',\tau]$, $k = 1,...,q$  (which are the resultants of the system (\ref{sup5.8})
with respects to the variables $w'$ and $\omega$) with the following property: for
every solution  
$(Z_0,z'_0,\zeta_0,\tau_0,w_0,\omega_0)$ of (\ref{sup5.8}), such that  
the leading coefficients of  polynomials in (\ref{sup5.8}) are different from zero at
$(Z_0,z'_0)$ and $(\zeta_0,\tau_0)$, there exists a neighborhood $W$ of
$(Z_0,z'_0,\zeta_0,\tau_0,w_0,\omega_0)$ such that for
$(Z,z',\zeta,\tau,w,\omega)\in W$  the system (\ref{sup5.8}) is equivalent to

$$\hat{Q}_j(Z,z',w'_j) = 0, \tilde{{Q}_j}(\zeta,\tau,\omega_j) = 0,
R_k(Z,z',\zeta,\tau) = 0$$

In order to construct $R_k$ explicitely,  consider the resultant
$R_k^1(Z,z',\tau,w'_2,...,w'_N,\omega)$ of $\hat{Q}_1(Z,z',w_1)$ and
$P'_k(z',\tau,w',\omega)$ with respect to the variable $w'_1$. Since $a \in (X
\cap U) \backslash \Sigma$, the leading coefficient of $\hat{Q}_1(Z,z',w_1)$ does not
vanish at
$(a,g(a))$. Hence, the replacing of $P'_k$ by $R_k^1$ in the system (\ref{sup5.8})
does not change its solutions in a neighborhood of $(a,f(a))$. Consider the resultant
$R_k^2(Z,z',\tau,w'_2,...,w'_N,\omega)$ of $R_k^1$ and $\hat{Q}_2$ with respect to
$w'_2$ , etc. This shows that every $P'_k$ in (\ref{sup5.8}) can be replaced by the
function $R_k^N(Z,z',\tau,\omega)$ which is a polynomial in $z',\tau,\omega$ with
coefficients holomorphic in $Z$ in a neighborhood of $p$. Consider
the the resultant $R_k^{N+1}(Z,z',\zeta,\tau,\omega_2,...,\omega_N)$ of
$\tilde{Q}_1$ and $R_k^N$ with respect to $\omega_1$ , etc. and repeat the same
arguments to eliminate variables $\omega_j$. The functions $R_k := R_k^{2N}$ satisfy
the desired properties. This completes the proof of proposition.

\bigskip

Now we are able to show (\ref{supsup5.1}). It follows from Proposition \ref{pro5.3}
that there exist  polynomials $R_k(Z,\overline{Z},z',\overline{z}') \in {\cal
O}^{\R}(X \cap U)[z',\overline{z}'], k = 1,...,q$ such that for every $a \in (X \cap
U) \backslash \Sigma$ there are  neighborhoods $V_1$ of $(a,g(a))$ and $V_2$ of
$h(a)$ such that 
$${\cal A}_{{}_p{\bf f}} \cap (X \times Y) \cap (V_1 \times V_2) = ({\cal A}_{{}_p{\bf
f}}\vert X) \cap \{ (Z,z') \in V_1: R_{k}(Z,\overline{Z}, z',
\overline{z}') = 0 , Z \in X, k = 1,...,q\}$$  
 Since the graph $\Gamma_{f}$ is contained in ${\cal A}_{{}_p{\bf f}} \cap (X
\times Y)$, we have $R_k(Z,\overline{Z},g,\overline{g}') = 0$ on
$(X \cap U) \backslash \Sigma$. The system $\{ {\bf g}_1,...,{\bf g}_m \}$ is 
algebraically independent  over ${\cal M}_p(X)$ and $\Sigma \subset (X \cap U)$ is a
${\cal P}$ - set, so Corollary \ref{cor4.5} implies that every $R_k$ represents
zero in ${\cal O}_p^{\R}(X)[z',\overline{z}']$. Therefore, for any $a \in X \backslash
\Sigma$, close enough to $p$ there exists a neighborhood
$V$ of the point $(a,f(a))$ in $\cc^n \times \cc^N$ such that $(({\cal
A}_{{}_p{\bf f}}\vert X) \cap V) \cap (X\times Y) = (({\cal A}_{{}_p{\bf f}}\vert X)
\cap V)$. Hence, $\pi'(({\cal A}_{{}_p{\bf f}}\vert X) \cap V)$ is contained in $Y$,
which completes the proof of the Main Theorem.

\bigskip

The estimate of the transcendence degree and the embedding of the graph of the
mapping in a complex analytic variety in Theorem \ref{theoB}  follow from the Main
Theorem and the remark in section 2. In particular, let the dimension of this variety
be equal to $n$ i.e. $tr.deg {}_p{\bf f} = 0$.  Since ${}_p{\bf f}$ admits a one-side
holomorphic extension and is $C^{\infty}$, the well known result of Bedford - Bell
\cite{BeB} implies that ${}_p{\bf f}$ is real analytic.

\bigskip

Bernard COUPET and Alexandre SUKHOV: 
 LATP, CNRS/ UMR n$^\circ$ 6632, CMI, Universit\'e de Provence, 39, rue Joliot Curie,
13453 Marseille cedex 13, France.

\bigskip 

Sergey PINCHUK: Department of Mathematics,
 Indiana University, Bloomington ,Indiana 47405, USA.

\end{document}